\documentclass[11pt]{article}
\usepackage[T1]{fontenc}
\usepackage{lmodern}
\usepackage[margin=27mm]{geometry}
\usepackage{amsmath,amssymb,amsthm,mathrsfs,mathtools,stmaryrd}
\usepackage{enumitem}
\usepackage{tikz-cd}
\usepackage[backend=biber,style=alphabetic,sorting=nyt,giveninits=true,
  maxnames=5,doi=true,isbn=false,url=false]{biblatex}
\usepackage{hyperref}
\hypersetup{colorlinks=true,linkcolor=blue,citecolor=blue,urlcolor=blue,
  pdftitle={Demazure--Lusztig type operators and correspondences on Kato Catalan varieties},
  pdfauthor={Taiyo Kodama}}
\allowdisplaybreaks[3]
\mathtoolsset{showonlyrefs=true}
\setlist[enumerate,1]{label=(\roman*),ref=\roman*}
\theoremstyle{plain}
\newtheorem{theorem}{Theorem}[section]
\newtheorem{proposition}[theorem]{Proposition}
\newtheorem{lemma}[theorem]{Lemma}
\newtheorem{corollary}[theorem]{Corollary}
\theoremstyle{definition}
\newtheorem{definition}[theorem]{Definition}

\theoremstyle{remark}
\newtheorem{remark}[theorem]{Remark}
\numberwithin{equation}{section}
\DeclareMathOperator{\Hom}{Hom}

\newcommand{\C}{\mathbb C}
\newcommand{\Z}{\mathbb Z}

\newcommand{\PP}{\mathbb P}
\newcommand{\id}{\mathrm{id}}
\newcommand{\BS}{\mathbf S}
\newcommand{\Hdeg}{\mathbf H_{\mathrm{deg}}}
\newcommand{\Gr}{\mathrm{Gr}}
\newcommand{\CO}{\mathcal O}
\newcommand{\CK}{\mathcal K}
\newcommand{\CB}{\mathcal B}
\newcommand{\CN}{\mathcal N}
\newcommand{\Fn}{\mathfrak n}
\newcommand{\pt}{\mathrm{pt}}
\newcommand{\SX}{\mathscr X}
\newcommand{\rfunc}{\mathtt r}
\newcommand{\hfunc}{\mathtt h}

\begin{document}
\title{Demazure--Lusztig type operators\\
and correspondences on Kato Catalan varieties}
\author{Taiyo Kodama\thanks{Graduate School of Science, Kyoto University. Email: \texttt{kodama.taiyo.43c@st.kyoto-u.ac.jp}.}}
\date{}
\maketitle
\begin{abstract}
We study the integral equivariant cohomology rings of Kato Catalan varieties $\SX_\Psi$ associated with arbitrary root ideals $\Psi$ in type~$A$.
These varieties are smooth projective compactifications of the vector subbundles $T^*_\Psi\CB$ of the cotangent bundle of the complete flag variety $\CB$.
Let $K$ denote the product of the general linear group (or its maximal torus) with the loop rotation torus $\C^\times$.
We give an explicit presentation of $H^*_K(\SX_\Psi;\Z)$ with relations expressed in terms of divided difference and Demazure--Lusztig type operators.
Furthermore, for every simple reflection preserving $\Psi\cup(-\Psi)$, we construct a geometric realization of the associated Demazure--Lusztig type operator on $H^*_K(\SX_\Psi;\Z)$ via correspondences.
These correspondences yield a faithful graded action of the subalgebra of the degenerate affine Hecke algebra generated by these reflections and the polynomial subalgebra.
For the maximal root ideal $\Psi=\Delta^+$, our construction recovers the cohomological counterpart of the $K$-theoretic Steinberg correspondences of Cautis and Kamnitzer.
\end{abstract}
\smallskip


\section{Introduction}\label{sec_intro}

Let $G=\mathrm{GL}_n(\C)$, let $B\subset G$ be the Borel subgroup of upper triangular matrices, let $T\subset B$ be the diagonal maximal torus, and write $\CB=G/B$ for the complete flag variety. Let $\mathfrak g=\mathfrak{gl}_n(\C)$ be the Lie algebra of $G$, and let $\Delta^+$ denote the set of positive roots of $(G,T)$ corresponding to $B$, with root spaces $\mathfrak g_\alpha\subset\mathfrak g$. Via the trace pairing on $\mathfrak g$, the cotangent bundle $T^*\CB$ is identified with $G\times^B\bigoplus_{\alpha\in\Delta^+}\mathfrak g_\alpha$. Associated with each root ideal $\Psi\subset\Delta^+$, Kato \cite[Theorems~3.10 and 4.4]{kato2023geometric} constructed a smooth projective compactification $\SX_\Psi$ of the $G$-equivariant vector subbundle
\[
 T^*_\Psi\CB=G\times^B\bigoplus_{\alpha\in\Psi}\mathfrak g_\alpha\subset T^*\CB.
\]
We refer to $\SX_\Psi$ as a \emph{Kato Catalan variety}.\footnote{In \cite{kato2023geometric}, Kato did not assign a name to $\SX_\Psi$; in private communication, he mentioned that he had originally intended to call it a \emph{Catalan variety}, but refrained from doing so as that terminology was already in use elsewhere.} These varieties were introduced to provide a geometric realization of Catalan functions and, via the endpoint morphism to the affine Grassmannian $\Gr$ and the geometric Satake equivalence, of chromatic quasisymmetric polynomials \cite{kato2023geometric,kato2024geometric,kamnitzer2024fixed}. From the viewpoint of $\Gr$, Kamnitzer described $\SX_\Psi$ as a reverse Hessenberg convolution variety and exhibited it as a tower of projective bundles over a point \cite[Proposition~2.1 and Remark~2.3]{kamnitzer2024fixed}.

For the maximal root ideal $\Psi=\Delta^+$, let $\CN\subset\mathfrak{gl}_n(\C)$ denote the nilpotent cone. The composite of the open embedding $T^*\CB\hookrightarrow\SX_{\Delta^+}$ with the endpoint morphism $\mathsf m_{\Delta^+}:\SX_{\Delta^+}\to\Gr$ is the composite of the Springer resolution $T^*\CB\to\CN$ with the classical lattice embedding $\CN\hookrightarrow\Gr$. This lattice construction originates with Lusztig \cite{lusztig1981green} and was developed further by Ng\^o \cite[Lemma~2.3.1]{ngo1999faisceaux} and Mirkovi\'c--Vybornov \cite[Theorem~3.1 and \S3.3]{mirkovic2003quiver}.
Cautis and Kamnitzer \cite{cautis2018quantum} studied the equivariant $K$-theory of $\SX_{\Delta^+}$, describing its $K$-ring via the convolution tower and constructing an affine Hecke algebra action from a smooth divisor and a lattice-forgetting fibration.
They compared the resulting correspondences with the Steinberg correspondences on the open cotangent bundle $T^*\CB$ \cite[\S11]{cautis2018quantum}.

In this paper, we study the integral equivariant cohomology rings of $\SX_\Psi$ for arbitrary root ideals $\Psi\subset\Delta^+$ and the operators defined by geometric correspondences, retaining the loop rotation action of $S=\C^\times$ throughout. We give a uniform operator presentation of $H^*_K(\SX_\Psi;\Z)$ for $K=T\times S$ and $G\times S$, and construct geometric realizations of the Demazure--Lusztig type operators for all simple reflections preserving $\Psi\cup(-\Psi)$.

To state our results precisely, fix $n\geq2$ and set $[n]=\{1,\dots,n\}$. Let $S=\C^\times$ act on $\C(\!(t)\!)^n$ by $v(t)\mapsto v(at)$ for $a\in S$. Root ideals $\Psi\subset\Delta^+$ are in bijection with \emph{reverse Hessenberg functions}, namely weakly increasing maps $\rfunc:[n]\to\{0,\dots,n-1\}$ with $\rfunc(k)<k$ for all $k\in[n]$, via
\[
 \Psi_{\rfunc}=\{\epsilon_a-\epsilon_b\in\Delta^+\mid a\leq\rfunc(b)\},
\]
where $\epsilon_a$ is the $a$-th coordinate character of $T$. In the lattice model, $\SX_{\rfunc}=\SX_{\Psi_{\rfunc}}$ parametrizes flags $L_0\subset L_1\subset\cdots\subset L_n$ of $\C\llbracket t\rrbracket$-lattices in $\C(\!(t)\!)^n$ satisfying
\[
 L_0=t\C\llbracket t\rrbracket^n,\qquad
 \dim_\C(L_k/L_{k-1})=1,\qquad
 tL_k\subset L_{\rfunc(k)}\quad(1\leq k\leq n).
\]
Set
\[
 \Lambda=\Lambda_{T\times S}=\Z[t_1,\dots,t_n,\epsilon],\qquad
 R=\Lambda[x_1,\dots,x_n],\qquad
 f(z)=\prod_{a=1}^n(z+t_a),
\]
where $t_1,\dots,t_n\in H^2_T(\pt;\Z)$ are the equivariant weights of the standard representation $\C^n$, $\epsilon=c_1^S(\C_1)$ for the one-dimensional representation $\C_1$ of $S$ of weight $1$, and $\deg t_a=\deg\epsilon=\deg x_k=2$. Let $\tau_i$ be the $\Lambda$-algebra involution of $R$ swapping $x_i$ and $x_{i+1}$, and define the divided difference $\partial_i$ and the Demazure--Lusztig type operator $\BS_i$ on $R$ by
\[
 \partial_i g=\frac{g-\tau_i g}{x_i-x_{i+1}},\qquad
 \BS_i=\tau_i-\epsilon\partial_i
       =\id+(x_{i+1}-x_i-\epsilon)\partial_i
\]
(see \cite[\S12]{ginzburg1998geometric}).

For $1\leq j\leq n$, let $c_j=c_j^G(\C^n)\in H_G^{2j}(\pt;\Z)$ denote the equivariant Chern classes of the standard representation of $G$, also viewed as classes for $G\times S$ by pullback. Under restriction to $T\times S$, these classes map to the elementary symmetric polynomials $e_j(t_1,\dots,t_n)$.

\begin{theorem}[Theorem~\ref{thm:equiv_coho_xdelta} and Corollary~\ref{cor:equiv_coho_xdelta_G}]\label{thm:intro_ring}
For any reverse Hessenberg function $\rfunc$, there is a canonical isomorphism of graded $\Lambda$-algebras
\[
 H^*_{T\times S}(\SX_{\rfunc};\Z)\cong R/I_{\rfunc},\qquad
 I_{\rfunc}=(D_1(f),\dots,D_n(f)),
\]
where, with empty operator products understood as the identity,
\[
 D_k(f)=\partial_{k-1}\cdots\partial_{\rfunc(k)+1}
             \BS_{\rfunc(k)}\cdots\BS_1 f(x_1).
\]
Replacing $\Lambda$ by $\Z[c_1,\dots,c_n,\epsilon]$ and $f(z)$ by $z^n+c_1z^{n-1}+\cdots+c_n$ yields the presentation of $H^*_{G\times S}(\SX_{\rfunc};\Z)$. Under these isomorphisms for $K=T\times S$ and $G\times S$, $x_k$ maps to $-c_1^K(L_k/L_{k-1})$.
\end{theorem}

The proof of Theorem~\ref{thm:intro_ring} uses a polynomial recurrence to identify the reversed polynomial associated with $D_k(f)$, modulo $(D_1(f),\dots,D_{k-1}(f))$, with the Chern polynomial of $t^{-1}L_{\rfunc(k)}/L_{k-1}$. Since $D_k(f)\in\Lambda[x_1,\dots,x_k]$ is monic of degree $d_k=n-k+1+\rfunc(k)$ in $x_k$, the relations are triangular; successive division yields normal forms in the monomial basis $\{x_1^{e_1}\cdots x_n^{e_n}\mid 0\leq e_k<d_k\}$. Equivariant integration then extracts the coefficient of $x_1^{d_1-1}\cdots x_n^{d_n-1}$ (Corollary~\ref{cor:normal_forms}).

Next, let $S'$ be the set of simple reflections $s_i=(i,i+1)$ in the Weyl group $\mathfrak S_n$ preserving $\Psi\cup(-\Psi)$, and let $\Hdeg'$ be the subalgebra of the degenerate affine Hecke algebra $\Hdeg$ generated by $S'$ and $\Z[x_1,\dots,x_n,\epsilon]$. For $K=T\times S$ or $G\times S$, write $\Phi_K:R_K=\Lambda_K[x_1,\dots,x_n]\twoheadrightarrow H_K^*(\SX_\Psi;\Z)$ for the presentation map.

\begin{theorem}[Theorem~\ref{thm:geometric_hecke} and Corollary~\ref{cor:deg_hecke_action_xdelta}]\label{thm:intro_hecke}
For each $s_i\in S'$, define the degree-zero endomorphism $\mathcal S_i$ of $H_K^*(\SX_\Psi;\Z)$ as follows:
\begin{enumerate}
\item If $\alpha_i=\epsilon_i-\epsilon_{i+1}\in\Psi$, let $j_i:Y_i\hookrightarrow\SX_\Psi$ be the smooth Cartier divisor defined by the condition $tL_{i+1}\subset L_{i-1}$, and let $\pi_i:Y_i\to B_i$ be the morphism to the corresponding partial lattice flag variety $B_i$ obtained by forgetting $L_i$. Then $\pi_i$ is an equivariant $\PP^1$-bundle, and we set
\[
 \mathcal S_i=\id+j_{i*}\pi_i^*\pi_{i*}j_i^*.
\]
\item If $\alpha_i\notin\Psi$, then forgetting $L_i$ defines a global equivariant $\PP^1$-bundle $\pi_i:\SX_\Psi\to B_i$, and we set
\[
 \mathcal S_i=\id+c_1^K(\Omega^1_{\SX_\Psi/B_i}\otimes\C_{-1})\cup\pi_i^*\pi_{i*},
\]
where $\C_{-1}$ is the one-dimensional representation on which $G$ acts trivially and $S$ acts with weight $-1$.
\end{enumerate}
In both cases, $\mathcal S_i\Phi_K(g)=\Phi_K(\BS_i g)$ for all $g\in R_K$. Consequently, $\ker\Phi_K$ is stable under $\BS_i$, and the operators $\mathcal S_i$ ($s_i\in S'$), together with multiplication by $x_j$ and $\epsilon$, define a graded representation of $\Hdeg'$ on $H_K^*(\SX_\Psi;\Z)$.
\end{theorem}

Restriction to the open subbundle $T^*_\Psi\CB$ defines a surjective $\Hdeg'$-module homomorphism $H_K^*(\SX_\Psi;\Z)\twoheadrightarrow H_K^*(T^*_\Psi\CB;\Z)$ sending each $x_i$ to the corresponding tautological class on $T^*_\Psi\CB$ (also denoted by $x_i$). For $G\times S$, the target is the polynomial representation $\Z[x_1,\dots,x_n,\epsilon]$. This proves faithfulness over $\Z[\epsilon]$ for $G\times S$; restriction to $T\times S$ gives faithfulness for that group as well (Corollary~\ref{cor:faithful_hecke}). Equivariant duality then gives an embedding of $\Hdeg'$ into the convolution algebra of correspondence classes on $\SX_\Psi\times\SX_\Psi$ (Corollary~\ref{thm:faithful_convolution_embedding}). For $\Psi=\Delta^+$, the fiber product $\overline Z_i=Y_i\times_{B_i}Y_i$ is the closure of the simple Steinberg component $Z_i\subset T^*\CB\times_{\CN}T^*\CB$, relating this action to the classical Steinberg correspondences (Subsection~\ref{sec:steinberg}).

\paragraph{Organization of the paper.}
Section~\ref{sec:deg_hecke} establishes the polynomial recurrences for $D_k(f)$, and Section~\ref{sec_xpsi} proves Theorem~\ref{thm:intro_ring}.
Section~\ref{sec:pushpull} computes the Gysin pushforwards for inclusions $\SX_{\rfunc}\hookrightarrow\SX_{\rfunc'}$ (Proposition~\ref{prop:general_gysin}) and the push--pull operators for lattice-forgetting morphisms, and Section~\ref{sec:geometric_hecke} proves Theorem~\ref{thm:intro_hecke}.
Section~\ref{sec_TX_incl} treats open restriction, faithfulness, and the comparison with Steinberg correspondences. Section~\ref{sec:examples} works out explicit examples in ranks two and three.

\section{Combinatorics and the degenerate affine Hecke algebra}\label{sec:deg_hecke}

\subsection{Root ideals and reverse Hessenberg functions}\label{ss:root_ideals_hess_funcs}

Let $n\geq2$ and set $[n]=\{1,\dots,n\}$. Let $G=\mathrm{GL}_n(\C)$, with Borel subgroup $B$
of upper triangular matrices and maximal torus $T$ of diagonal matrices.
Write $\mathfrak g=\mathfrak{gl}_n(\C)$ and $\mathfrak t=\operatorname{Lie}(T)$
for the Lie algebras of $G$ and $T$, respectively.
Let $\epsilon_i\in\mathfrak t^*$ be the $i$-th coordinate character of $T$.
The root system of $(G,T)$ and its positive roots are
\[
  \Delta=\{\epsilon_i-\epsilon_j\mid i\neq j\},
  \qquad
  \Delta^+=\{\epsilon_i-\epsilon_j\mid i<j\}.
\]
A subset $\Psi\subset\Delta^+$ is a \emph{root ideal} if
$(\Psi+\Delta^+)\cap\Delta^+\subset\Psi$.
Any root ideal $\Psi$ is uniquely encoded by a
\emph{reverse Hessenberg function}, that is, a weakly increasing map
$\rfunc:[n]\to\{0,\dots,n-1\}$ satisfying $\rfunc(j)<j$ for all $j\in[n]$.
The bijection between root ideals and reverse Hessenberg functions is given by
\[
  \Psi_{\rfunc}
  =\{\epsilon_i-\epsilon_j\in\Delta^+\mid i\leq\rfunc(j)\},
  \qquad
  \rfunc_\Psi(j)
  =\max\bigl(\{i\mid\epsilon_i-\epsilon_j\in\Psi\}\cup\{0\}\bigr).
\]
The Weyl group $W=N_G(T)/T$ is identified with the symmetric group $\mathfrak S_n$,
where the simple reflection $s_i=(i,i+1)$ interchanges $i$ and $i+1$.

\begin{remark}[Hessenberg functions and reverse Hessenberg functions]\label{rem:hessenberg_comparison}
In combinatorics and the geometry of Hessenberg varieties, it is more common to encode a root ideal $\Psi\subset\Delta^+$ by a \emph{Hessenberg function}, namely a weakly increasing map $\hfunc:[n]\to[n]$ satisfying $\hfunc(i)\geq i$ for all $i\in[n]$. A Hessenberg function $\hfunc$ determines a root ideal by $\Psi_\hfunc=\{\epsilon_i-\epsilon_j\in\Delta^+\mid j>\hfunc(i)\}$. The bijection between Hessenberg functions $\hfunc$ and reverse Hessenberg functions $\rfunc$ corresponding to the same root ideal is given by the equivalence $i\leq\rfunc(j)\iff j>\hfunc(i)$, or explicitly by
\[
  \rfunc(j)=\max\bigl(\{i\in[n]\mid \hfunc(i)<j\}\cup\{0\}\bigr),
  \qquad
  \hfunc(i)=\max\{j\in[n]\mid\rfunc(j)<i\}.
\]
Although Hessenberg functions are more standard in the combinatorics literature, we follow Kamnitzer \cite[\S2.1]{kamnitzer2024fixed} in using reverse Hessenberg functions throughout this paper, as they allow us to write the defining lattice condition of $\SX_{\rfunc}$ and the rank of each stage of the projective bundle tower directly in terms of $\rfunc(j)$ (see Section~\ref{sec_xpsi}).
\end{remark}

\subsection{The degenerate affine Hecke algebra and polynomial operators}\label{ss:deg_hecke_algs}

Let $\epsilon$ be a formal parameter of degree two.
We recall the degenerate affine Hecke algebra (cf.\ Ginzburg \cite[\S 12]{ginzburg1998geometric}).

\begin{definition}\label{def:deg-hecke}
The \emph{degenerate affine Hecke algebra} $\Hdeg$ is the associative algebra over $\Z[\epsilon]$
generated by pairwise commuting elements $x_1,\dots,x_n$ and simple reflections $s_1,\dots,s_{n-1}$
satisfying the Coxeter relations of $\mathfrak S_n$ together with the cross relations
\[
  s_jx_i-x_{s_j(i)}s_j=
  \begin{cases}
    -\epsilon & \text{if } i=j,\\
    \epsilon & \text{if } i=j+1,\\
    0 & \text{if } i\neq j,j+1,
  \end{cases}
  \qquad (1\leq i\leq n,\ 1\leq j<n).
\]
The algebra $\Hdeg$ is graded by $\deg s_j=0$ and $\deg x_i=\deg\epsilon=2$.
\end{definition}

Let $\Lambda$ be a commutative $\Z[\epsilon]$-algebra and set $R=\Lambda[x_1,\dots,x_n]$.
Let $\tau_i$ denote the $\Lambda$-algebra automorphism of $R$ interchanging
$x_i$ and $x_{i+1}$ and fixing all other variables.
We define the divided difference operator $\partial_i$ and the
Demazure--Lusztig type operator $\BS_i$ on $R$ by
\[
  \partial_i f=\frac{f-\tau_i f}{x_i-x_{i+1}},
  \qquad
  \BS_i f=\tau_i f-\epsilon\partial_i f.
\]
The divided difference operators satisfy the nil-Hecke relations and the
twisted Leibniz rule
\begin{equation}\label{eq:divided_difference_leibniz}
  \partial_i(fg)=(\partial_i f)g+(\tau_i f)\partial_i g.
\end{equation}
The operators $\BS_i$ satisfy the Coxeter relations of $\mathfrak S_n$, and together
with multiplication by $x_i$, they define the polynomial representation
\begin{equation}\label{eq:deg-hecke-action}
  s_i\longmapsto\BS_i,\qquad
  x_i\longmapsto x_i
\end{equation}
of $\Hdeg$ on $R$ (see \cite[Proposition~12.2]{ginzburg1998geometric}).
Note that $\BS_i$ can be rewritten in the form
\begin{equation}\label{eq:BS-expansion}
 \BS_i=\id+(x_{i+1}-x_i-\epsilon)\partial_i.
\end{equation}
When $R$ is graded with $\deg\epsilon=2$ and $\deg x_i=2$, the
operator $\partial_i$ has degree $-2$ and $\BS_i$ has degree zero.

\begin{lemma}\label{lem:polynomial_faithful}
The polynomial representation of $\Hdeg$ on
$\Z[x_1,\dots,x_n,\epsilon]$ is faithful.
\end{lemma}

\begin{proof}
For $w\in W$, let $s_w$ be the product of simple generators along a reduced expression.
By the cross relations, every element of $\Hdeg$ can be written as
$\sum_{w\in W}p_w s_w$ with $p_w\in\Z[x_1,\dots,x_n,\epsilon]$.
If this element acts by zero on $\Z[x_1,\dots,x_n,\epsilon]$, then its image in
$\C\otimes_{\Z}\Hdeg$ acts by zero on $\C[x_1,\dots,x_n,\epsilon]$ and is therefore zero
by \cite[Proposition~12.2(ii)]{ginzburg1998geometric}.
The complex PBW decomposition \cite[Definition~12.1(i)]{ginzburg1998geometric}
then gives $p_w=0$ for every $w\in W$, so the element is zero in $\Hdeg$.
\end{proof}

\subsection{Polynomials for the cohomology presentation}\label{ss:lemmas_polynomials}

We now establish polynomial identities that express the relations arising from the
projective bundle formula for the tower defining $\SX_{\rfunc}$ (see Section~\ref{sec_xpsi}).

Let $f(z)\in\Lambda[z]$ be a monic polynomial of degree $n$.
For $0\leq j<k\leq n$, define the polynomial $P_{k,j}(f)\in R$ by
\[
  P_{k,j}(f)
  =\partial_{k-1}\cdots\partial_{j+1}
       \BS_j\cdots\BS_1 f(x_1),
\]
where empty products of operators are understood to be the identity.
For a reverse Hessenberg function $\rfunc$, set
\[
  D_k(f)=P_{k,\rfunc(k)}(f),\qquad
  I_{\rfunc}^{[k]}(f)=(D_1(f),\dots,D_k(f))\subset R,
  \qquad I_{\rfunc}^{[0]}(f)=0,
\]
and write $I_{\rfunc}(f)=I_{\rfunc}^{[n]}(f)$.

\begin{lemma}\label{lem:polynomial_recurrence}
The polynomial $P_{k,j}(f)$ belongs to $\Lambda[x_1,\dots,x_k]$ and is monic
of degree $n-k+1+j$ in $x_k$.
Moreover, for all $1\leq j<k$, it satisfies the recurrence
\begin{equation}\label{eq:triangular_polynomial_recurrence}
  P_{k,j}(f)
  =(x_k-x_j-\epsilon)P_{k,j-1}(f)+P_{k-1,j-1}(f).
\end{equation}
\end{lemma}

\begin{proof}
If $g\in\Lambda[x_1,\dots,x_i]$ is monic of degree $d>0$ in $x_i$,
then $\BS_i g$ and $\partial_i g$ are monic in $x_{i+1}$ of degrees $d$ and $d-1$,
respectively.
In the definition of $P_{k,j}(f)$, we apply $j$ Demazure--Lusztig type operators to $f(x_1)$,
followed by $k-j-1$ divided differences.
Thus $P_{k,j}(f)\in\Lambda[x_1,\dots,x_k]$ is monic of degree $n-k+1+j$ in $x_k$.

For fixed $j$, the case $k=j+1$ of
\eqref{eq:triangular_polynomial_recurrence} follows from
\eqref{eq:BS-expansion} applied to $P_{j,j-1}(f)$.
For $k>j+1$, we apply $\partial_{k-1}$ to the recurrence for $k-1$.
Since $P_{k-2,j-1}(f)\in\Lambda[x_1,\dots,x_{k-2}]$ is independent of $x_{k-1}$ and $x_k$,
it is annihilated by $\partial_{k-1}$, and the Leibniz rule
\eqref{eq:divided_difference_leibniz} gives the result, using
$\partial_{k-1}(x_{k-1}-x_j-\epsilon)=1$ and
$\tau_{k-1}(x_{k-1}-x_j-\epsilon)=x_k-x_j-\epsilon$.
\end{proof}

We introduce the reversed polynomials in an auxiliary variable $u$:
\[
  F(u)=u^nf(u^{-1}),\qquad
  d_{k,j}=n-k+1+j,\qquad
  \widehat P_{k,j}(u)
    =u^{d_{k,j}}P_{k,j}(f)\big|_{x_k=u^{-1}}.
\]
By Lemma~\ref{lem:polynomial_recurrence}, $\widehat P_{k,j}(u)\in \Lambda[x_1,\dots,x_{k-1}][u]$
is a polynomial in $u$ of degree at most $d_{k,j}$ with constant term $1$.
The following lemma computes $\widehat P_{k,j}(u)$ modulo $I_{\rfunc}^{[k-1]}(f)$.

\begin{lemma}\label{lem:operator_chern_relation}
Let $1\leq k\leq n$ and $\rfunc(k)\leq j<k$.
\begin{enumerate}
  \item\label{item:Pk}
    $P_{k,j}(f)\in I_{\rfunc}^{[k]}(f)$.
  \item\label{item:Qk}
    In $(R/I_{\rfunc}^{[k-1]}(f))\llbracket u\rrbracket$,
    \begin{equation}\label{eq:operator_chern_series}
      \widehat P_{k,j}(u)
      =
      F(u)\,
      \frac{\prod_{i=1}^{j}(1-(x_i+\epsilon)u)}
           {\prod_{i=1}^{k-1}(1-x_i u)}.
    \end{equation}
\end{enumerate}
\end{lemma}

\begin{proof}
We prove (i) by induction on $k$. The base case $k=1$ is
$P_{1,0}(f)=D_1(f)\in I_{\rfunc}^{[1]}(f)$.
For $k\geq 2$, the boundary case $j=\rfunc(k)$ holds by definition:
$P_{k,\rfunc(k)}(f)=D_k(f)\in I_{\rfunc}^{[k]}(f)$.
For $j>\rfunc(k)$, the monotonicity of $\rfunc$ implies $j-1\geq\rfunc(k)\geq\rfunc(k-1)$,
so $P_{k-1,j-1}(f)\in I_{\rfunc}^{[k-1]}(f)$ by induction on $k$.
The recurrence \eqref{eq:triangular_polynomial_recurrence} then gives
$P_{k,j}(f)\in I_{\rfunc}^{[k]}(f)$ by induction on $j$.

We prove (ii) by induction on $k$. The base case $k=1$ is $\widehat P_{1,0}(u)=F(u)$.
For $k\geq 2$, we distinguish two cases.

If $j\leq k-2$, then $P_{k,j}(f)=\partial_{k-1}P_{k-1,j}(f)$.
Since $P_{k-1,j}(f)\in\Lambda[x_1,\dots,x_{k-1}]$ is independent of $x_k$, applying $\tau_{k-1}$ substitutes $x_k$ for $x_{k-1}$, so
\[
  u^{d_{k-1,j}}(\tau_{k-1}P_{k-1,j}(f))\big|_{x_k=u^{-1}}=\widehat P_{k-1,j}(u).
\]
Set $d=d_{k,j}$, so that $d_{k-1,j}=d+1$.
Expanding the divided difference and substituting $x_k=u^{-1}$ gives
\[
  (1-x_{k-1}u)\widehat P_{k,j}(u)
  =\widehat P_{k-1,j}(u)-u^{d+1}P_{k-1,j}(f).
\]
Since $j\geq\rfunc(k)\geq\rfunc(k-1)$, the second term vanishes modulo
$I_{\rfunc}^{[k-1]}(f)$ by part~(i).
Applying the induction hypothesis to $\widehat P_{k-1,j}(u)$ and dividing by
$1-x_{k-1}u$ yields \eqref{eq:operator_chern_series}.

If $j=k-1$, then $P_{k,k-1}(f)=\BS_{k-1}P_{k-1,k-2}(f)$.
Since $d_{k,k-1}=d_{k-1,k-2}=n$, applying the same substitution to
$P_{k-1,k-2}(f)\in\Lambda[x_1,\dots,x_{k-1}]$ gives
$u^n(\tau_{k-1}P_{k-1,k-2}(f))\big|_{x_k=u^{-1}}=\widehat P_{k-1,k-2}(u)$.
Using \eqref{eq:BS-expansion} and setting $x_k=u^{-1}$, we find
\[
  (1-x_{k-1}u)\widehat P_{k,k-1}(u)
  =(1-(x_{k-1}+\epsilon)u)\widehat P_{k-1,k-2}(u)
  +\epsilon u^{n+1}P_{k-1,k-2}(f).
\]
Since $k-2\geq\rfunc(k-1)$ by the condition $\rfunc(k-1)<k-1$, the second term
again vanishes modulo $I_{\rfunc}^{[k-1]}(f)$ by part~(i).
Applying the induction hypothesis to $\widehat P_{k-1,k-2}(u)$ and dividing by
$1-x_{k-1}u$ gives \eqref{eq:operator_chern_series}.
\end{proof}

In Subsection~\ref{ss:equiv_coho_xpsi}, we identify $\widehat P_{k,\rfunc(k)}(u)$,
modulo the preceding relations, with the Chern polynomial at the $k$-th stage of the projective bundle tower.

\subsection{Simple-reflection symmetries of root ideals}

The geometric realization of the operators $\BS_i$ in Section~\ref{sec:geometric_hecke}
requires the reflection $s_m$ to preserve the set $\Psi\cup(-\Psi)$.
In terms of the reverse Hessenberg function, this symmetry is characterized as follows.

\begin{lemma}\label{lem:symmetry_of_root_ideal}
Let $\Psi=\Psi_{\rfunc}$ and $1\leq m<n$.
The set $\Psi\cup(-\Psi)$ is $s_m$-invariant if and only if
$\rfunc(k)\neq m$ for all $k\geq m+2$ and one of the following two
mutually exclusive conditions holds:
\begin{enumerate}
  \item $\rfunc(m)=m-1$ and $\rfunc(m+1)=m$;\label{item:case1A}
  \item $\rfunc(m)=\rfunc(m+1)$.\label{item:case1B}
\end{enumerate}
In case~\ref{item:case1A}, $\epsilon_m-\epsilon_{m+1}\in\Psi$;
in case~\ref{item:case1B}, $\epsilon_m-\epsilon_{m+1}\notin\Psi$.
\end{lemma}

\begin{proof}
Aside from the simple root $\epsilon_m-\epsilon_{m+1}$, whose sign is reversed,
$s_m$ interchanges the roots in each of the pairs
\[
  \{\epsilon_i-\epsilon_m,\epsilon_i-\epsilon_{m+1}\}\quad(i<m),
  \qquad
  \{\epsilon_m-\epsilon_k,\epsilon_{m+1}-\epsilon_k\}\quad(k>m+1).
\]
For the first family, membership in $\Psi$ agrees within each pair
if and only if $\min(\rfunc(m+1),m-1)=\rfunc(m)$.
Since $\rfunc(m)\leq\rfunc(m+1)\leq m$, this equality is equivalent to
either (i) $\rfunc(m+1)=m$ and $\rfunc(m)=m-1$, or (ii) $\rfunc(m+1)=\rfunc(m)\le m-1$.
For the second family, membership agrees within each pair if and only if $\rfunc(k)\neq m$.
All remaining positive roots are fixed by $s_m$.
Since $\Psi\cup(-\Psi)$ is symmetric under negation, this proves the claim.
\end{proof}

\section{Equivariant cohomology rings of \texorpdfstring{$\SX_\Psi$}{X\_Psi}}\label{sec_xpsi}
\subsection{Conventions for equivariant cohomology}\label{sec_pre_equ_coho}

Unless specified otherwise, all cohomology groups are taken with integral coefficients
in the sense of Borel equivariant cohomology.
Let $G = \mathrm{GL}_n(\C)$, let $T \subset G$ be the maximal torus of diagonal matrices,
and let $S = \C^\times$ act on the one-dimensional representation $\C_1 \cong \C$ by
the standard character of weight $1$ (loop rotation).
For a linear algebraic group $K$, we write $\Lambda_K = H_K^*(\pt; \Z)$.
Our primary coefficient rings are
\begin{equation}\label{eq:coefficient_rings}
 \Lambda_{T\times S} = \Z[t_1,\dots,t_n,\epsilon], \qquad
 \Lambda_{G\times S} = \Z[c_1,\dots,c_n,\epsilon],
\end{equation}
where $t_a\in H^2_T(\pt;\Z)$ is the equivariant first Chern class of the $a$-th coordinate line of $\C^n$ (corresponding to the character $\epsilon_a$), $\epsilon = c_1^S(\C_1)$, and $c_j = c_j^G(\C^n)$,
with gradings $\deg t_a = \deg\epsilon = 2$ and $\deg c_j = 2j$.
The restriction homomorphism $\Lambda_{G\times S}\hookrightarrow\Lambda_{T\times S}$ identifies
$c_j$ with the $j$-th elementary symmetric polynomial $e_j(t_1,\dots,t_n)$.
For an equivariant vector bundle $E$, we denote its equivariant Chern polynomial by
$c_u^K(E) = \sum_{j \ge 0} c_j^K(E) u^j$; this extends by multiplicativity to virtual bundles
as a formal power series in $u$.

\begin{lemma}\label{lem:borel}
Let $A$ be a complex algebraic torus (possibly trivial), and let $X$ be a
smooth complex quasi-projective variety equipped with an algebraic
$G\times A$-action.
Assume that the ordinary cohomology $H^*(X;\Z)$ is free over $\Z$ and
concentrated in even degrees.
Then the restriction homomorphism induces isomorphisms of graded algebras
\[
 H^*_{G\times A}(X;\Z)
 \otimes_{\Lambda_{G\times A}}\Lambda_{T\times A}
 \xrightarrow{\sim}H^*_{T\times A}(X;\Z),
 \qquad
 H^*_{G\times A}(X;\Z)
 \xrightarrow{\sim}H^*_{T\times A}(X;\Z)^W.
\]
Here the Weyl group $W = \mathfrak S_n$ acts on $H^*_{T\times A}(X;\Z)$ via the normalizer $N_G(T)$,
permuting the classes $t_a$, fixing $\Lambda_A$, and fixing all classes in the image of $H^*_{G\times A}(X;\Z)$.
\end{lemma}

\begin{proof}
Set $K = G \times A$ and $K_T = T \times A$.
Since $H^*(X;\Z)$ and $H^*(BK;\Z) \cong \Lambda_K$ are free $\Z$-modules
concentrated in even degrees, the Serre spectral sequence for the Borel fibration
$X \to EK \times_K X \to BK$ collapses at the $E_2$-page.
By the Leray--Hirsch theorem, any homogeneous integral basis $\{b_j\}$ of $H^*(X;\Z)$
lifts to a free $\Lambda_K$-basis $\{\widetilde b_j\}$ of $H_K^*(X;\Z)$.
These lifts restrict to a $\Lambda_{K_T}$-basis of $H_{K_T}^*(X;\Z)$,
yielding the base-change isomorphism.

Since the action of $W$ on $EK \times_{K_T} X$ covers the identity map on $EK \times_K X$,
each restricted basis element $\operatorname{res}(\widetilde b_j)$ is $W$-invariant.
Expanding an arbitrary class $\alpha \in H_{K_T}^*(X;\Z)$ in this basis, we see that $\alpha$
is $W$-invariant if and only if its coefficients lie in
$\Lambda_{K_T}^W = \operatorname{res}(\Lambda_K)$.
This identifies $H_K^*(X;\Z)$ with the invariant subalgebra $H_{K_T}^*(X;\Z)^W$.
\end{proof}

We adopt the convention that the projective bundle $\PP(E)$ parametrizes one-dimensional
linear subspaces (lines) in a vector bundle $E$ of rank $r$ over $X$.
Setting $x = -c_1^K(\CO_{\PP(E)}(-1))$, the equivariant projective bundle formula gives
\begin{equation}\label{eq:projective_bundle_formula}
 H_K^*(\PP(E);\Z) \cong
 H_K^*(X;\Z)[x]\Big/\Big(\sum_{j=0}^r c_j^K(E)x^{r-j}\Big).
\end{equation}
The classes $1, x, \dots, x^{r-1}$ form a free basis of $H_K^*(\PP(E);\Z)$ as an
$H_K^*(X;\Z)$-module, and the fiber integration map $\pi_*:H_K^*(\PP(E);\Z)\to H_K^{*-2(r-1)}(X;\Z)$
is normalized by $\pi_*(x^{r-1}) = 1$ and $\pi_*(x^i) = 0$ for $0\leq i < r-1$
(cf.\ \cite[Proposition~2.6.1]{anderson2024equivariant}).

\subsection{Lattice flags}\label{ss:def_xpsi}

Let $\CO=\C\llbracket t\rrbracket$ and $\CK=\C(\!(t)\!)$.
A \emph{lattice} in $\CK^n$ is a free $\CO$-submodule of rank $n$.
The set of all lattices in $\CK^n$ is parametrized by the \emph{affine Grassmannian}
$\Gr=G(\CK)/G(\CO)$ of $G=\mathrm{GL}_n$.

Let $S=\C^\times$. The semidirect product $G(\CO)\rtimes S$ acts on $\CK^n$ by
\[
  (g,a)\cdot v(t) = g \cdot v(at)
\]
for $g\in G(\CO)$, $a\in S$, and $v(t)\in \CK^n$, inducing an action of $G(\CO)\rtimes S$ on $\Gr$.

Following Kato \cite{kato2023geometric} and Kamnitzer \cite{kamnitzer2024fixed},
the \emph{Kato Catalan variety} $\SX_{\rfunc}$ associated with a reverse
Hessenberg function $\rfunc$ is defined by
\[
 \SX_{\rfunc}=\left\{L_\bullet=(L_0,\dots,L_n)\ \middle|\
 \begin{gathered}
 t\CO^n=L_0\subset L_1\subset\cdots\subset L_n,\quad L_i\in\Gr,\\
 \dim_\C(L_i/L_{i-1})=1,\quad tL_i\subset L_{\rfunc(i)}
 \quad(1\leq i\leq n)
 \end{gathered}\right\}.
\]
We also write $\SX_\Psi$ for $\SX_{\rfunc}$ to emphasize the corresponding root ideal $\Psi$.
The variety $\SX_{\rfunc}$ is stable under the action of $G(\CO)\rtimes S$.
The endpoint morphism is
\[
  \mathsf m_\Psi:\SX_\Psi\longrightarrow\Gr,\qquad
  (L_0,\dots,L_n)\longmapsto t^{-1}L_n.
\]
This morphism is $G(\CO)\rtimes S$-equivariant.

\begin{remark}[The initial lattice]\label{rem:lattice_normalization}
We choose $L_0=t\CO^n$ so that the quotient
$t^{-1}L_0/L_0\cong\C^n$ carries the trivial $S$-action (trivial loop rotation).
Consequently, the tautological classes $x_i=-c_1^K(L_i/L_{i-1})$ restrict to the standard
dual flag classes on the open vector subbundle $T^*_\Psi\CB$ (see Section~\ref{sec_TX_incl}).
The alternative normalization $\widetilde L_0=\CO^n$ used in
\cite[\S2.1]{kamnitzer2024fixed} and \cite[Lemma~6.1]{cautis2018quantum}
is obtained by setting $\widetilde L_i=t^{-1}L_i$; the defining lattice inclusions
remain unchanged. Under this shift, the endpoint morphism is simply given by
$\widetilde L_\bullet\mapsto\widetilde L_n=\mathsf m_\Psi(L_\bullet)$, and the tautological lines are twisted
by the $S$-representation $\C_{-1}$ of weight $-1$:
\[
 \widetilde L_i/\widetilde L_{i-1}
 \cong(L_i/L_{i-1})\otimes\C_{-1},\qquad
 \widetilde x_i=-c_1^K(\widetilde L_i/\widetilde L_{i-1})
   =x_i+\epsilon.
\]
Accordingly, in the $\CO^n$ convention, the polynomial $f_K(z)$ below is replaced by $f_K(z-\epsilon)$.
\end{remark}

We define a partial order on the set of reverse Hessenberg functions by setting
$\rfunc \leq \rfunc'$ if $\rfunc(i) \leq \rfunc'(i)$ for all $i \in [n]$.
For $\rfunc \leq \rfunc'$, the defining inclusions yield a natural
closed embedding $\SX_{\rfunc} \hookrightarrow \SX_{\rfunc'}$.

\subsection{A presentation of the equivariant cohomology ring}\label{ss:equiv_coho_xpsi}

For $K=T\times S$ or $G\times S$, set $R_K=\Lambda_K[x_1,\dots,x_n]$ and define the polynomial
\[
  f_K(z)=
  \begin{cases}
    \prod_{a=1}^n(z+t_a) & \text{for } K=T\times S,\\
    z^n+c_1z^{n-1}+\cdots+c_n & \text{for } K=G\times S.
  \end{cases}
\]
Using the operators introduced in Subsection~\ref{ss:deg_hecke_algs}, we define the ideal
\[
  I_{\rfunc,K}=(D_k(f_K)\mid1\leq k\leq n)\subset R_K.
\]
We now identify these algebraic relations with those given by the projective
bundle formula at each stage of the tower.

\begin{theorem}\label{thm:equiv_coho_xdelta}
Let $\rfunc$ be a reverse Hessenberg function and set
$f(z)=\prod_{a=1}^n(z+t_a)\in\Lambda_{T\times S}[z]$.
Then there is an isomorphism of graded $\Lambda_{T\times S}$-algebras
\[
  H^*_{T\times S}(\SX_{\rfunc};\Z) \cong \Lambda_{T\times S}[x_1,\dots,x_n]/I_{\rfunc,T\times S},
\]
where $I_{\rfunc,T\times S}=(D_k(f)\mid1\leq k\leq n)$ and the
polynomials $D_k(f)$ are defined in Subsection~\ref{ss:lemmas_polynomials}.
Under this isomorphism, $x_i$ maps to $-c_1^{T\times S}(L_i/L_{i-1})$.
\end{theorem}

\begin{proof}
Write $\Lambda=\Lambda_{T\times S}$ and suppress the equivariant superscript in the notation for Chern classes.
For each $0\leq k\leq n$, let $\SX_{\rfunc}^{[k]}$ denote the variety of partial flags $(L_0,\dots,L_k)$
satisfying the defining conditions of $\SX_{\rfunc}$ up to index $k$.
Forgetting the terminal lattice $L_k$ yields a tower of morphisms
\[
  \SX_{\rfunc}=\SX_{\rfunc}^{[n]}\longrightarrow\SX_{\rfunc}^{[n-1]}\longrightarrow\cdots
  \longrightarrow\SX_{\rfunc}^{[0]}=\pt.
\]
By \cite[Proposition~2.1]{kamnitzer2024fixed}, the $k$-th stage is
$\SX_{\rfunc}^{[k]}\cong\PP(\mathscr V_k)$, where, writing $a=\rfunc(k)$,
\[
 \mathscr V_k=t^{-1}L_a/L_{k-1},\qquad \operatorname{rk}\mathscr V_k=d_k=n-k+1+a.
\]
Its tautological line bundle is $L_k/L_{k-1}$. We also write $x_k$ for $-c_1(L_k/L_{k-1})$.
We proceed by induction on $k$ to prove that
\[
  H^*_{T\times S}(\SX_{\rfunc}^{[k]};\Z)
  \cong\Lambda[x_1,\dots,x_k]/(D_1(f),\dots,D_k(f)).
\]
The base case $k=0$ is immediate since $\SX_{\rfunc}^{[0]}=\pt$.
Assuming that the presentation holds for $\SX_{\rfunc}^{[k-1]}$, we determine the Chern
relation introduced by $\PP(\mathscr V_k)$.
From the two short exact sequences
\[
  0\longrightarrow L_{k-1}/L_0\longrightarrow t^{-1}L_a/L_0\longrightarrow\mathscr V_k\longrightarrow0,
\]
\[
  0\longrightarrow t^{-1}L_0/L_0\longrightarrow t^{-1}L_a/L_0
    \longrightarrow t^{-1}L_a/t^{-1}L_0\longrightarrow0,
\]
we obtain the following identity in equivariant $K$-theory:
\[
  [\mathscr V_k]
   =[\C^n]+\sum_{i=1}^{a}[t^{-1}L_i/t^{-1}L_{i-1}]
            -\sum_{i=1}^{k-1}[L_i/L_{i-1}].
\]
Here $\C^n=t^{-1}L_0/L_0$ is the standard representation of $T$ with trivial
$S$-action, having equivariant weights $t_1,\dots,t_n$.
Multiplication by $t^{-1}$ twists each tautological line by the $S$-representation $\C_{-1}$ of weight $-1$ (whose equivariant first Chern class is $-\epsilon$),
so $c_1(t^{-1}L_i/t^{-1}L_{i-1})=-(x_i+\epsilon)$.
The Whitney sum formula therefore gives
\begin{equation}\label{eq:chern_polynomial_Vk}
  c_u(\mathscr V_k)
  =\prod_{b=1}^n(1+t_bu)\,
    \frac{\prod_{i=1}^{a}(1-(x_i+\epsilon)u)}
         {\prod_{i=1}^{k-1}(1-x_i u)}
\end{equation}
in $H^*_{T\times S}(\SX_{\rfunc}^{[k-1]};\Z)\llbracket u\rrbracket$.

Since $F(u)=\prod_{b=1}^n(1+t_bu)$, Lemma~\ref{lem:operator_chern_relation}(ii)
with $j=a=\rfunc(k)$, together with the induction hypothesis, gives
\[
  \sum_{j=0}^{d_k}c_j(\mathscr V_k)u^j
   =\widehat P_{k,a}(u)
   =u^{d_k}D_k(f)\big|_{x_k=u^{-1}}.
\]
Comparing coefficients of $u$ yields
\[
  D_k(f)=\sum_{m=0}^{d_k}c_m(\mathscr V_k)x_k^{d_k-m}
\]
over $H^*_{T\times S}(\SX_{\rfunc}^{[k-1]};\Z)$.
By the equivariant projective bundle formula \eqref{eq:projective_bundle_formula},
the relation imposed at stage $k$ is precisely $D_k(f)=0$, completing the induction.
\end{proof}

The vector bundle sequences in the proof above also underlie the $K$-theoretic
presentation for the maximal function $\rfunc(k)=k-1$ in
\cite[Lemmas~8.13--8.14]{cautis2018quantum}.

Since the ordinary projective bundle formula shows that $H^*(\SX_{\rfunc};\Z)$
is a free $\Z$-module concentrated in even degrees, Lemma~\ref{lem:borel} applies.

\begin{corollary}\label{cor:equiv_coho_xdelta_G}
There are canonical isomorphisms
\[
  H^*_{G\times S}(\SX_{\rfunc};\Z)
  \cong H^*_{T\times S}(\SX_{\rfunc};\Z)^W
  \cong R_{G\times S}/I_{\rfunc,G\times S},
\]
with $x_i=-c_1^{G\times S}(L_i/L_{i-1})$.
Here the normalizer action of the Weyl group $W$ permutes the $t_i$ and fixes all
$x_i$ and $\epsilon$.
\end{corollary}

\begin{proof}
The first isomorphism follows from Lemma~\ref{lem:borel}.
For the second, the projective bundle tower in the proof of
Theorem~\ref{thm:equiv_coho_xdelta} is $G\times S$-equivariant.
Writing $a=\rfunc(k)$, the same exact sequences give
\[
 c_u^{G\times S}(\mathscr V_k)
 =(1+c_1u+\cdots+c_nu^n)
 \frac{\prod_{i=1}^{a}(1-(x_i+\epsilon)u)}
      {\prod_{i=1}^{k-1}(1-x_i u)}.
\]
Assuming the presentation through stage $k-1$,
Lemma~\ref{lem:operator_chern_relation}(ii) gives
\[
 c_u^{G\times S}(\mathscr V_k)
 =u^{d_k}D_k(f_{G\times S})\big|_{x_k=u^{-1}},
 \qquad d_k=n-k+1+a,
\]
since $u^nf_{G\times S}(u^{-1})=1+c_1u+\cdots+c_nu^n$.
Comparing coefficients identifies $D_k(f_{G\times S})$ with the projective
bundle relation. Starting from $\SX_{\rfunc}^{[0]}=\pt$, we therefore obtain
the presentation by induction on $k$ using the equivariant projective bundle formula.
The stated Weyl group action follows from Lemma~\ref{lem:borel} and the
$G\times S$-equivariance of the tautological lines.
\end{proof}

\begin{corollary}\label{cor:coho_xdelta_R}
With integer coefficients, restriction to a subtorus $T'\subset T\times S$ (possibly trivial) is given by base change:
\[
 H^*_{T'}(\SX_{\rfunc};\Z)
 \cong H^*_{T\times S}(\SX_{\rfunc};\Z)
       \otimes_{\Lambda_{T\times S}}\Lambda_{T'}.
\]
In particular, setting $f(z)=\prod_{a=1}^n(z+t_a)$, one has
\[
 H_T^*(\SX_{\rfunc};\Z)
 \cong\frac{\Z[t_1,\dots,t_n,x_1,\dots,x_n]}
 {\bigl(\partial_{k-1}\cdots\partial_{\rfunc(k)+1}
                      f(x_{\rfunc(k)+1})\mid1\leq k\leq n\bigr)}.
\]
The ordinary cohomology ring $H^*(\SX_{\rfunc};\Z)$ is obtained from this $T$-equivariant presentation by setting all $t_a=0$.
For $K=T\times S$ or $G\times S$ and every commutative ring $A$, one also has
\[
 H_K^*(\SX_{\rfunc};A)
 \cong H_K^*(\SX_{\rfunc};\Z)\otimes_\Z A.
\]
\end{corollary}

\begin{proof}
For base change, choose homogeneous equivariant lifts of an integral basis of ordinary cohomology.
Under restriction to $T'$, these classes continue to restrict to the same ordinary basis;
applying the Leray--Hirsch theorem to both Borel fibrations identifies the restriction with the stated tensor product.
Setting $\epsilon=0$ replaces each operator $\BS_i$ by the transposition $\tau_i$,
and $\tau_j\cdots\tau_1f(x_1)=f(x_{j+1})$, giving the displayed relations.

The coefficient-change isomorphism follows from the universal coefficient
theorem, since the integral equivariant cohomology is free in each degree.
\end{proof}

\subsection{Monomial bases and equivariant integration}

For $\Psi=\Psi_{\rfunc}$, set
\[
 d_k=n-k+1+\rfunc(k),\qquad
 d=\sum_{k=1}^n(d_k-1)=\binom n2+|\Psi|.
\]

\begin{corollary}\label{cor:normal_forms}
For $K=T\times S$ or $G\times S$, the monomials
\[
 x_1^{e_1}\cdots x_n^{e_n},\qquad 0\leq e_k<d_k,
\]
form a $\Lambda_K$-basis of $H^*_K(\SX_{\rfunc};\Z)$.
Equivariant integration extracts the coefficient of
$x_1^{d_1-1}\cdots x_n^{d_n-1}$ in the expansion in this basis.
\end{corollary}

\begin{proof}
At the $k$-th stage of the projective bundle tower, the classes
$1,x_k,\dots,x_k^{d_k-1}$ form a basis over the cohomology of the preceding
stage. Iterating gives the stated $\Lambda_K$-basis.
Since integration along the $k$-th fiber sends $x_k^{d_k-1}$ to $1$,
\[
 \int_{\SX_{\rfunc}}x_1^{d_1-1}\cdots x_n^{d_n-1}=1.
\]
Every other basis monomial has degree less than $2d$ and hence has zero
integral. The formula follows by $\Lambda_K$-linearity.
\end{proof}

\subsection{Duality and correspondence classes}\label{ss:duality}

Let $X=\SX_{\rfunc}$, $d=\dim_\C X$, and $M_K=H_K^*(X;\Z)$ for $K=T\times S$ or $G\times S$.
The equivariant integration pairing
\[
  M_K\otimes_{\Lambda_K}M_K\longrightarrow\Lambda_K,
  \qquad (\alpha,\beta)\longmapsto\int_X\alpha\cup\beta,
\]
is perfect over $\Lambda_K$.
Indeed, choose homogeneous equivariant lifts $\widetilde a_i,\widetilde b_i$
of Poincar\'e-dual bases of $H^*(X;\Z)$, ordered by nondecreasing $\deg a_i$.
By the Leray--Hirsch theorem, these lifts are $\Lambda_K$-bases of $M_K$. Since
\[
 \deg\int_X\widetilde a_i\widetilde b_j=\deg a_i-\deg a_j,
\]
the pairing matrix is lower triangular with diagonal entries $1$:
negative-degree entries vanish, and degree-zero entries are the ordinary intersection numbers.
External products of these lifts likewise lift an integral basis of $H^*(X\times X;\Z)$,
so the Leray--Hirsch theorem gives the canonical K\"unneth isomorphism
\[
  M_K\otimes_{\Lambda_K}M_K
  \xrightarrow{\sim}H_K^*(X\times X;\Z).
\]

Let $\operatorname{pr}_1,\operatorname{pr}_2:X\times X\to X$ denote the projections.
A correspondence class $\gamma\in H_K^{2d+r}(X\times X;\Z)$ defines a $\Lambda_K$-linear
endomorphism of degree $r$ on $M_K$ by
\begin{equation}\label{eq:correspondence_action}
 \mathcal T_\gamma(u)=(\operatorname{pr}_1)_*
       (\gamma\cup\operatorname{pr}_2^*u).
\end{equation}
Any closed $K$-invariant subvariety $Z\subset X\times X$ of codimension $c$ determines
a correspondence class $[Z]\in H_K^{2c}(X\times X;\Z)$ via its equivariant fundamental cycle
and Poincar\'e duality on $X\times X$.
For the projection morphisms $p_{ab}:X\times X\times X\to X\times X$, the
convolution product is
\begin{equation}\label{eq:convolution_product}
 \gamma*\eta=(p_{13})_*(p_{12}^*\gamma\cup p_{23}^*\eta).
\end{equation}
By the projection formula, $\mathcal T_{\gamma*\eta}=\mathcal T_\gamma\circ\mathcal T_\eta$,
and the diagonal class $[\Delta_X]$ acts as the identity.
Set $\mathcal C_X^r=H_K^{r+2d}(X\times X;\Z)$, graded so that convolution has degree zero.
The correspondence action induces a canonical isomorphism of graded algebras
\begin{equation}\label{eq:ambient_endomorphisms}
 \mathcal C_X^\bullet\xrightarrow{\sim}
 \operatorname{End}_{\Lambda_K}^\bullet(M_K).
\end{equation}
Indeed, under the K\"unneth isomorphism it sends $\alpha\otimes\beta$ to
$u\mapsto\alpha\int_X\beta u$, which is an isomorphism by the perfectness of the pairing.

\section{Gysin maps and push--pull operators}\label{sec:pushpull}

The cohomology ring presentation of Subsection~\ref{ss:equiv_coho_xpsi} provides
an explicit polynomial description of the Gysin maps for lattice inclusions
and the push--pull operators for forgetful morphisms.
Throughout this section, we work integrally over the coefficient ring $\Lambda_K$,
with $K=T\times S$ or $K=G\times S$, and retain the notation
$R_K=\Lambda_K[x_1,\dots,x_n]$, $f_K$, and $I_{\rfunc,K}$, with $x_i=-c_1^K(L_i/L_{i-1})$.

For $\rfunc'\leq\rfunc$, the tautological line bundles on $\SX_{\rfunc}$ restrict along the inclusion $j:\SX_{\rfunc'}\hookrightarrow\SX_{\rfunc}$ to those on $\SX_{\rfunc'}$. Consequently, $I_{\rfunc,K}\subset I_{\rfunc',K}$ and the pullback $j^*$ is identified with the natural quotient map $j^*:R_K/I_{\rfunc,K}\twoheadrightarrow R_K/I_{\rfunc',K}$.

\subsection{Gysin maps for divisor embeddings}

Let $\rfunc$ and $\rfunc'$ be reverse Hessenberg functions such that $\rfunc'$ is obtained
from $\rfunc$ by decrementing a single value:
\[
  \rfunc'(m)=\rfunc(m)-1 \quad\text{for some } m\in\{2,\dots,n\},
  \qquad \rfunc'(k)=\rfunc(k) \quad(k\neq m).
\]
Set $a=\rfunc(m)$ and let $i:\SX_{\rfunc'}\hookrightarrow\SX_{\rfunc}$
be the canonical inclusion.

\begin{proposition}\label{prop:inclusion_gysin}
The closed embedding $i:\SX_{\rfunc'}\hookrightarrow\SX_{\rfunc}$ realizes
$\SX_{\rfunc'}$ as a smooth Cartier divisor in $\SX_{\rfunc}$, with normal bundle
isomorphic to the restriction of the equivariant line bundle
\[
  \mathcal E=\Hom(tL_m/tL_{m-1},L_a/L_{a-1}).
\]
Under the presentations of Theorem~\ref{thm:equiv_coho_xdelta} and Corollary~\ref{cor:equiv_coho_xdelta_G}, the Gysin pushforward map
$i_*:H_K^q(\SX_{\rfunc'};\Z)\to H_K^{q+2}(\SX_{\rfunc};\Z)$ is given by
\[
  i_*[f]=[(x_m-x_a-\epsilon)f] \qquad (f\in R_K).
\]
\end{proposition}

\begin{proof}
Weak monotonicity of $\rfunc'$ gives $\rfunc(m-1)=\rfunc'(m-1)\leq\rfunc'(m)=a-1$,
so $tL_{m-1}\subset L_{a-1}$ holds everywhere on $\SX_{\rfunc}$.
Together with $tL_m\subset L_a$, the lattice inclusion induces an equivariant homomorphism
$tL_m/tL_{m-1}\to L_a/L_{a-1}$, defining a section $s\in\Gamma(\SX_{\rfunc},\mathcal E)$
whose zero locus is precisely $\SX_{\rfunc'}$.
At stage $m$ of the projective bundle tower, $tL_{m-1}\subset L_{a-1}$
makes $t^{-1}L_{a-1}/L_{m-1}$ a corank-one subbundle of
$t^{-1}L_a/L_{m-1}$, and $s=0$ defines the smooth relative hyperplane
\[
 \PP(t^{-1}L_{a-1}/L_{m-1})\subset\PP(t^{-1}L_a/L_{m-1}).
\]
Let $p_m:\SX_{\rfunc}\to\SX_{\rfunc}^{[m]}$ forget the lattices after $L_m$.
Since only the $m$-th lattice condition changes, $\SX_{\rfunc'}$ is the
scheme-theoretic inverse image of this hyperplane under $p_m$.
The map $p_m$ is a composition of projective bundle projections and is therefore smooth.
Thus $s$ is transverse to the zero section and the normal bundle is $i^*\mathcal E$.

Since multiplication by $t$ twists by $\C_1$, the first Chern classes of the two quotient lines are $c_1^K(tL_m/tL_{m-1})=-x_m+\epsilon$ and $c_1^K(L_a/L_{a-1})=-x_a$, giving
$c_1^K(\mathcal E)=x_m-x_a-\epsilon$.
Because $i_*(1)=c_1^K(\mathcal E)$, the projection formula yields $i_*(i^*\alpha)=\alpha\cup c_1^K(\mathcal E)=(x_m-x_a-\epsilon)\alpha$
for all $\alpha\in H_K^*(\SX_{\rfunc};\Z)$.
Since $i^*$ is surjective, every class on $\SX_{\rfunc'}$ is of the form $i^*[f]=[f]$,
yielding $i_*[f]=[(x_m-x_a-\epsilon)f]$.
\end{proof}

\subsection{Forgetting a lattice}

We next consider the forgetful morphism that drops the intermediate lattice $L_m$ from a flag.
Fix an index $1\leq m<n$, and assume that
\begin{equation}\label{eq:pushpull_hypotheses}
  \rfunc(m)=\rfunc(m+1),\qquad
  \rfunc(k)\neq m\quad(k\geq m+2).
\end{equation}
Let $K_m=\{0,\dots,n\}\setminus\{m\}$.
By \eqref{eq:pushpull_hypotheses}, $\rfunc(k)\neq m$ for every $k$,
so all targets $L_{\rfunc(k)}$ belong to the retained flag.
We may therefore define the partial flag variety
\[
  \SX_{\rfunc,m}^P
  =\left\{(L_k)_{k\in K_m}\ \middle|\;
  \begin{gathered}
    L_k\in\Gr,\quad L_0=t\CO^n,\quad L_k\subset L_\ell\ (k<\ell),\\
    \dim_\C(L_k/L_0)=k,\quad
    tL_k\subset L_{\rfunc(k)}\ (k\in K_m\setminus\{0\})
  \end{gathered}\right\},
\]
and denote by $\pi_m:\SX_{\rfunc}\to\SX_{\rfunc,m}^P$ the morphism forgetting $L_m$.

\begin{proposition}\label{prop:pushpull_projective_bundle}
Under \eqref{eq:pushpull_hypotheses}, the morphism $\pi_m$ realizes $\SX_{\rfunc}$ as the
equivariant projective bundle $\PP(\mathcal V_m)$ over $\SX_{\rfunc,m}^P$, where
$\mathcal V_m=L_{m+1}/L_{m-1}$ is a rank-two equivariant vector bundle.
In particular, the base variety $\SX_{\rfunc,m}^P$ is smooth and projective.
\end{proposition}

\begin{proof}
Set $a=\rfunc(m)=\rfunc(m+1)\leq m-1$. For any partial flag in $\SX_{\rfunc,m}^P$,
we have $tL_{m+1}\subset L_a\subset L_{m-1}$, so every intermediate lattice
$L_{m-1}\subset L_m\subset L_{m+1}$ automatically satisfies $tL_m\subset tL_{m+1}\subset L_a$.
Since no other condition refers to $L_m$ by \eqref{eq:pushpull_hypotheses},
a lattice $L_m$ completes the partial flag if and only if $L_{m-1}\subset L_m\subset L_{m+1}$ with $\dim_\C(L_m/L_{m-1})=1$.
Such lattices correspond bijectively to lines in $L_{m+1}/L_{m-1}$,
identifying $\SX_{\rfunc}\cong\PP(\mathcal V_m)$ equivariantly.
The base is smooth because the projective bundle is locally a product with $\PP^1$
and $\SX_{\rfunc}$ is smooth.
For projectivity, the $n$-dimensional $\CO$-module $L_n/L_0$ is annihilated by $t^n$,
so every retained subspace $L_k/L_0$ lies in the fixed vector space $t^{-n}L_0/L_0$.
The incidence and lattice conditions are closed conditions in the corresponding
partial flag variety, so $\SX_{\rfunc,m}^P$ is projective.
\end{proof}

\subsection{The geometric divided difference operator}

We now relate the push--pull operator $\pi_m^*\pi_{m*}$ to the divided difference
operator $\partial_m$ on $R_K$.

\begin{theorem}\label{thm:pushpull_divided_difference}
Assume \eqref{eq:pushpull_hypotheses}.
Let $\Phi_K:R_K\twoheadrightarrow H_K^*(\SX_{\rfunc};\Z)$ be the presentation homomorphism.
Then for every $f\in R_K$,
\[
  \pi_m^*\pi_{m*}\Phi_K(f)=\Phi_K(\partial_m f).
\]
Consequently, $I_{\rfunc,K}$ is stable under $\partial_m$ and $\tau_m$.
\end{theorem}

\begin{proof}
On $\SX_{\rfunc}$, the pulled-back bundle $\pi_m^*\mathcal V_m=L_{m+1}/L_{m-1}$ fits into
the tautological exact sequence $0\to L_m/L_{m-1}\to\pi_m^*\mathcal V_m\to L_{m+1}/L_m\to 0$.
Taking equivariant Chern classes yields $c_u^K(\pi_m^*\mathcal V_m)=(1-x_m u)(1-x_{m+1}u)$, so
\[
  \Phi_K(x_m+x_{m+1})=-\pi_m^*c_1^K(\mathcal V_m),\qquad
  \Phi_K(x_m x_{m+1})=\pi_m^*c_2^K(\mathcal V_m).
\]
For every index $k\notin\{m,m+1\}$, the line bundle $L_k/L_{k-1}$ is pulled back from $\SX_{\rfunc,m}^P$,
so $\Phi_K(x_k)\in\pi_m^*H_K^2(\SX_{\rfunc,m}^P;\Z)$.

Let $R_K^{s_m}\subset R_K$ be the subring of polynomials symmetric in $x_m,x_{m+1}$.
Every $f\in R_K$ can be written uniquely as $f=a+x_m b$ with $a,b\in R_K^{s_m}$,
whence $\partial_m f=b$ and both $\Phi_K(a)$ and $\Phi_K(b)$ lie in $\pi_m^*H_K^*(\SX_{\rfunc,m}^P;\Z)$.
Since $\pi_{m*}1=0$ and $\pi_{m*}x_m=1$, the projection formula gives
$\pi_m^*\pi_{m*}\Phi_K(f)=\Phi_K(b)=\Phi_K(\partial_m f)$.

Since $I_{\rfunc,K}=\ker\Phi_K$, if $f\in I_{\rfunc,K}$ then
$\Phi_K(\partial_m f)=\pi_m^*\pi_{m*}(0)=0$, so $\partial_m f\in I_{\rfunc,K}$.
From $\tau_m=\id-(x_m-x_{m+1})\partial_m$, it follows that $I_{\rfunc,K}$ is also stable under $\tau_m$.
\end{proof}

\subsection{Gysin maps for arbitrary inclusions}

We now generalize the Gysin formula of Proposition~\ref{prop:inclusion_gysin} to arbitrary inclusions.
Let $\rfunc'\leq\rfunc$ be reverse Hessenberg functions, and let
$j:\SX_{\rfunc'}\hookrightarrow\SX_{\rfunc}$ be the canonical closed embedding of codimension
$c=\sum_{k=1}^n(\rfunc(k)-\rfunc'(k))$.
Define the polynomial
\[
  E_{\rfunc',\rfunc}
  =\prod_{k=1}^n\prod_{a=\rfunc'(k)+1}^{\rfunc(k)}
       (x_k-x_a-\epsilon)\ \in R_K,
\]
which is homogeneous of cohomological degree $2c$.

\begin{proposition}\label{prop:general_gysin}
The Gysin pushforward $j_*:H_K^q(\SX_{\rfunc'};\Z)\to H_K^{q+2c}(\SX_{\rfunc};\Z)$ is given by
\[
  j_*[f]=[E_{\rfunc',\rfunc}f] \qquad (f\in R_K).
\]
Moreover, $j_*$ is injective.
\end{proposition}

\begin{proof}
We transition from $\rfunc$ to $\rfunc'$ by decreasing the value at index $k$ from $\rfunc(k)$ down to $\rfunc'(k)$ one unit at a time, proceeding in increasing order $k=1,\dots,n$.
At step $k$, the values at indices $h<k$ have already been reduced to $\rfunc'(h)$,
while those at $h>k$ remain $\rfunc(h)$.
Any intermediate value $v$ with $\rfunc'(k)\leq v\leq\rfunc(k)$ satisfies $0\leq v<k$ and
\[
 \begin{aligned}
  \rfunc'(k-1)&\leq\rfunc'(k)\leq v &&\text{if }k>1,\\
  v&\leq\rfunc(k)\leq\rfunc(k+1) &&\text{if }k<n.
 \end{aligned}
\]
Thus each intermediate function remains a reverse Hessenberg function.
This factors $j$ into $c$ inclusions of smooth Cartier divisors.
By Proposition~\ref{prop:inclusion_gysin}, each decrement from $a$ to $a-1$
at index $k$ contributes the factor $x_k-x_a-\epsilon$.
Functoriality of Gysin maps therefore gives $j_*[f]=[E_{\rfunc',\rfunc}f]$.

If $j_*\alpha=0$, the projection formula gives
\[
 \int_{\SX_{\rfunc'}}\alpha\cup j^*\beta
 =\int_{\SX_{\rfunc}}j_*\alpha\cup\beta=0
 \qquad\bigl(\beta\in H_K^*(\SX_{\rfunc};\Z)\bigr).
\]
Since $j^*$ is surjective and the integration pairing on $\SX_{\rfunc'}$
is perfect, we conclude that $\alpha=0$.
\end{proof}

\section{Geometric realizations of Demazure--Lusztig type operators}
\label{sec:geometric_hecke}

Let $\Psi=\Psi_{\rfunc}$ be a root ideal, and set $X=\SX_{\rfunc}$.
Let $S'$ be the set of simple reflections in the Weyl group $W$ preserving
the set $\Psi\cup(-\Psi)$, and let $\Hdeg'\subset\Hdeg$ denote the
subalgebra generated by $S'$ and the polynomial subalgebra $\Z[x_1,\dots,x_n,\epsilon]$.
By Lemma~\ref{lem:symmetry_of_root_ideal}, each $s_i\in S'$ satisfies
$\rfunc(k)\neq i$ for $k\geq i+2$. We use the following two constructions
according to whether $\epsilon_i-\epsilon_{i+1}$ belongs to $\Psi$.
Let $K=T\times S$ or $G\times S$.

\subsection{The two geometric constructions}

Fix $s_i\in S'$. We first consider the case where $\rfunc(i)=i-1$ and $\rfunc(i+1)=i$, so that the
simple root $\epsilon_i-\epsilon_{i+1}$ belongs to $\Psi$.
Define a reverse Hessenberg function $\rfunc_i^-$ by setting $\rfunc_i^-(i+1)=i-1$
and $\rfunc_i^-(k)=\rfunc(k)$ for all $k\neq i+1$.
Then $\rfunc_i^-$ is weakly increasing and determines morphisms
\[
  j_i:Y_i=\SX_{\rfunc_i^-}\hookrightarrow X,
  \qquad
  \pi_i:Y_i\longrightarrow B_i=\SX_{\rfunc_i^-,i}^P.
\]
By Proposition~\ref{prop:inclusion_gysin} (applied with $m=i+1$ and $a=i$), the embedding $j_i$
realizes $Y_i$ as a smooth Cartier divisor in $X$ whose normal bundle has first Chern class
$c_1^K(N_{Y_i/X})=j_i^*(x_{i+1}-x_i-\epsilon)$.
Moreover, since $\rfunc_i^-(i)=\rfunc_i^-(i+1)=i-1$ and $\rfunc_i^-(k)=\rfunc(k)\neq i$ for $k\geq i+2$,
the function $\rfunc_i^-$ satisfies \eqref{eq:pushpull_hypotheses} with $m=i$.
Thus Proposition~\ref{prop:pushpull_projective_bundle} identifies $\pi_i$ with the equivariant
projective bundle $\PP(L_{i+1}/L_{i-1})$ over the smooth base $B_i$.
Applying $j_i^*$, $\pi_{i*}$, $\pi_i^*$, and $j_{i*}$ in that order and adding the identity,
we define the degree-zero endomorphism
\begin{equation}\label{eq:hecke_case_A}
  \mathcal S_i=\id+j_{i*}\pi_i^*\pi_{i*}j_i^*
  \qquad\text{on }H^*_K(X;\Z).
\end{equation}
For the maximal root ideal $\Psi=\Delta^+$, this divisor and forgetful correspondence coincide
with the geometric setup of Cautis and Kamnitzer \cite[Lemma~9.3]{cautis2018quantum} in equivariant $K$-theory.

Next, for our fixed $s_i\in S'$, suppose that $\rfunc(i)=\rfunc(i+1)$, which by Lemma~\ref{lem:symmetry_of_root_ideal} is equivalent to $\epsilon_i-\epsilon_{i+1}\notin\Psi$.
In this case, the condition $\rfunc(k)\neq i$ for $k\geq i+2$ implies that the function $\rfunc$ itself
satisfies \eqref{eq:pushpull_hypotheses} with $m=i$.
Hence Proposition~\ref{prop:pushpull_projective_bundle} yields a projective line bundle on the entire variety:
\[
  \pi_i:X\longrightarrow B_i=\SX_{\rfunc,i}^P.
\]
Let $\mathcal E_i=\Omega^1_{X/B_i}\otimes\C_{-1}$ denote the relative cotangent line bundle
twisted by the one-dimensional representation where $G$ acts trivially and $S$ acts with weight $-1$.
Since the relative tangent bundle is $T_{X/B_i}\cong\Hom(L_i/L_{i-1},L_{i+1}/L_i)$, its dual twisted by $\C_{-1}$ has
equivariant Euler class
\[
  e_K(\mathcal E_i)=c_1^K(\mathcal E_i)
  =c_1^K(L_i/L_{i-1})-c_1^K(L_{i+1}/L_i)-\epsilon
  =x_{i+1}-x_i-\epsilon.
\]
Multiplying the push--pull operator by this degree-two class, we define the degree-zero endomorphism
\begin{equation}\label{eq:hecke_case_B}
  \mathcal S_i=\id+e_K(\mathcal E_i)\cup\pi_i^*\pi_{i*}
  \qquad\text{on }H^*_K(X;\Z).
\end{equation}

\subsection{Identification with the polynomial operators}

Let $\Phi_K:R_K\twoheadrightarrow H_K^*(X;\Z)$ denote the presentation homomorphism.

\begin{theorem}\label{thm:geometric_hecke}
For every $s_i\in S'$ and every polynomial $f\in R_K$, the geometric operators
\eqref{eq:hecke_case_A} and \eqref{eq:hecke_case_B} satisfy
\[
  \mathcal S_i\Phi_K(f)=\Phi_K(\BS_i f).
\]
In particular, the relation ideal $I_{\rfunc,K}$ is stable under $\BS_i$.
\end{theorem}

\begin{proof}
When $\epsilon_i-\epsilon_{i+1}\in\Psi$, Theorem~\ref{thm:pushpull_divided_difference},
applied to $Y_i$, identifies $\pi_i^*\pi_{i*}$ with $\partial_i$ on polynomial representatives.
Proposition~\ref{prop:inclusion_gysin}
then identifies $j_{i*}$ with multiplication by $x_{i+1}-x_i-\epsilon$.
When $\epsilon_i-\epsilon_{i+1}\notin\Psi$, the same push--pull theorem
applies on $X$, and the Euler class is $x_{i+1}-x_i-\epsilon$.
Thus in either case
\[
 (\mathcal S_i-\id)\Phi_K(f)
 =\Phi_K\bigl((x_{i+1}-x_i-\epsilon)\partial_i f\bigr)
 =\Phi_K\bigl((\BS_i-\id)f\bigr),
\]
by \eqref{eq:BS-expansion}. Stability of $I_{\rfunc,K}$ follows from $I_{\rfunc,K}=\ker\Phi_K$.
\end{proof}

\begin{corollary}\label{cor:deg_hecke_action_xdelta}
The assignments $s_i\mapsto\mathcal S_i$ for $s_i\in S'$, together with multiplication by
$x_j$ and $\epsilon$, define a graded representation of $\Hdeg'$ on $H^*_K(\SX_{\rfunc};\Z)$.
These representations extend compatibly to any commutative coefficient ring.
In particular, when $\Psi=\Delta^+$, they define an action of the full degenerate affine Hecke algebra $\Hdeg$.
\end{corollary}

\begin{proof}
The relations of the polynomial representation descend through $\Phi_K$
by Theorem~\ref{thm:geometric_hecke}. Coefficient extension follows from
Corollary~\ref{cor:coho_xdelta_R}.
\end{proof}

\subsection{Correspondences realizing the operators}

Let $d=\dim_\C X$, and define the closed subvariety
\[
  Z_{\rfunc}=\{(L_\bullet,L'_\bullet)\in X\times X\mid L_n=L'_n\},
\]
endowed with its reduced subscheme structure.
Because the forgetful morphism $\pi_i$ preserves the terminal lattice $L_n$ (as $i<n$),
all correspondences constructed below are supported on $Z_{\rfunc}$.

When $\epsilon_i-\epsilon_{i+1}\in\Psi$, let $C_i=Y_i\times_{B_i}Y_i$, regarded as a closed subvariety
of $X\times X$ via $j_i\times j_i$.
Since $\dim_\C Y_i=d-1$ and $\dim_\C B_i=d-2$, the fiber product $C_i$ has dimension $d$,
defining an equivariant fundamental class $[C_i]\in H_K^{2d}(X\times X;\Z)$.
We define the correspondence class
\[
  \Gamma_i=[\Delta_X]+[C_i]\in H_K^{2d}(X\times X;\Z).
\]
Let $q_1,q_2:C_i\to Y_i$ denote the two projections.
The fiber square
\[
\begin{tikzcd}
  C_i\arrow[r,"q_2"]\arrow[d,"q_1"'] & Y_i\arrow[d,"\pi_i"] \\
  Y_i\arrow[r,"\pi_i"'] & B_i
\end{tikzcd}
\]
is Cartesian with $\pi_i$ smooth and proper.
By flat base change for fiber integration, we have $q_{1*}q_2^*=\pi_i^*\pi_{i*}$.
Under the correspondence action \eqref{eq:correspondence_action}, the class $[C_i]$ therefore acts by
\[
  j_{i*}q_{1*}q_2^*j_i^*=j_{i*}\pi_i^*\pi_{i*}j_i^*,
\]
so $[\Delta_X]+[C_i]$ acts by \eqref{eq:hecke_case_A}.

When $\epsilon_i-\epsilon_{i+1}\notin\Psi$, the forgetful morphism $\pi_i:X\to B_i$ is defined on the entire variety $X$, and the corresponding
fiber product $C_i=X\times_{B_i}X$ has dimension $d+1$, so its fundamental class $[C_i]\in H_K^{2d-2}(X\times X;\Z)$ has
cohomological degree $2d-2$.
To obtain a degree-zero correspondence, we take the cup product with $\operatorname{pr}_1^*e_K(\mathcal E_i)$:
\[
  \Gamma_i=[\Delta_X]+\operatorname{pr}_1^*e_K(\mathcal E_i)\cup[C_i]\in H_K^{2d}(X\times X;\Z).
\]
By flat base change for $C_i=X\times_{B_i}X$ and the projection formula for $\operatorname{pr}_1$, the action of $\Gamma_i$ coincides with
\eqref{eq:hecke_case_B}.

\section{Restriction to the open subbundle and faithfulness}\label{sec_TX_incl}

\subsection{The open embedding}\label{ss:inclusion_TX_Xdelta_overview}

Let $\CB=G/B$ denote the variety of complete flags
$F_\bullet=(0=F_0\subset F_1\subset\cdots\subset F_n=\C^n)$ in $\C^n$, and let
$\CN\subset\mathfrak{gl}_n(\C)$ be the nilpotent cone.
For a root ideal $\Psi=\Psi_{\rfunc}$, set
\[
  \Fn(\Psi)=\bigoplus_{\alpha\in\Psi}\mathfrak g_\alpha,
  \qquad T^*_\Psi\CB=G\times^B\Fn(\Psi).
\]
Via the standard trace pairing identifying $T^*\CB$ with $G\times^B\Fn(\Delta^+)$,
the vector bundle $T^*_\Psi\CB$ embeds into $T^*\CB$ as the subbundle
\[
  T^*_\Psi\CB
  =\{(F_\bullet,\xi)\in\CB\times\CN
        \mid \xi F_i\subset F_{\rfunc(i)}\text{ for all }1\leq i\leq n\}.
\]
In particular, $T^*_{\Delta^+}\CB=T^*\CB$.
The group $G\times S$ acts on $T^*_\Psi\CB$ by
\[
  (g,a)\cdot(F_\bullet,\xi)=(gF_\bullet,a^{-1}g\xi g^{-1}),
\]
where loop rotation scales the cotangent fibers with weight $-1$.

We use the classical lattice chart of
\cite{lusztig1981green,ngo1999faisceaux,mirkovic2003quiver}, as described for
$T^*\CB$ in \cite[Lemma~11.2]{cautis2018quantum}, and restrict it to the
subbundle $T^*_\Psi\CB$ inside Kato's compactification \cite{kato2023geometric}.

Recall from Section~\ref{sec_xpsi} that points of $\SX_\Psi$ are lattice flags $L_\bullet$ with
$L_0=t\CO^n$.
Set $V=(\CK/t\CO)^n$. Each flag $L_\bullet$ is uniquely determined by the flag of
finite-dimensional subspaces $M_i=L_i/L_0\subset V$.
Because $tM_i\subset M_{\rfunc(i)}\subset M_{i-1}$, each $M_i$ is contained in $t^{1-n}\CO^n/t\CO^n$.
With respect to the decomposition $V=(\CK/t\CO)^n=\bigoplus_{h\geq 0}t^{-h}\C^n$, let $p:V\to\C^n$ be the projection onto the degree-zero component $\CO^n/t\CO^n\cong\C^n$.
For any $\xi\in\CN$, define the linear map
\[
  q_\xi:\C^n\longrightarrow V,\qquad
  q_\xi(v)=\sum_{h=0}^{n-1}t^{-h}\xi^h v.
\]
Because $\xi^n=0$, multiplication by $t$ on $V$ satisfies $t q_\xi=q_\xi\xi$, while $p\circ q_\xi=\id_{\C^n}$.

\begin{proposition}\label{prop:T^*X_incl}
Under the identification of each lattice flag $L_\bullet\in\SX_\Psi$
with its quotient flag $M_\bullet=(L_i/L_0)_{i=0}^n$ in $V$, the map
\[
  \iota_\Psi:T^*_\Psi\CB\longrightarrow\SX_\Psi,\qquad
  (F_\bullet,\xi)\longmapsto\bigl(q_\xi(F_i)\bigr)_{i=0}^n
\]
is a $G\times S$-equivariant open immersion.
Its image consists of those lattice flags for which the restriction
$p|_{M_n}:M_n\to\C^n$ is an isomorphism.
\end{proposition}

\begin{proof}
The identities $tq_\xi=q_\xi\xi$ and $p\circ q_\xi=\id_{\C^n}$ give
$tq_\xi(F_i)\subset q_\xi(F_{\rfunc(i)})$.
On the stated open locus, the inverse morphism is given by
\[
 F_i=p(M_i),\qquad \xi=(p|_{M_n})\circ(t|_{M_n})\circ(p|_{M_n})^{-1},
\]
where $t|_{M_n}:M_n\to M_n$ denotes multiplication by $t$.
Indeed, $t^nM_n=0$ and every $v\in M_n$ satisfies
$v=\sum_{h=0}^{n-1}t^{-h}p(t^hv)=q_\xi(p(v))$.
In particular, $M_i=q_\xi(F_i)$ for all $i$, and $\xi F_i=p(tM_i)\subset p(M_{\rfunc(i)})=F_{\rfunc(i)}$, so $(F_\bullet,\xi)\in T^*_\Psi\CB$.
Both assignments are algebraic and $G$-equivariant.
The identity $a\cdot q_\xi(v)=q_{a^{-1}\xi}(v)$ for $a\in S$ proves $S$-equivariance.
\end{proof}

\subsection{Restriction of the operators}
\label{ss:inclusion_TX_Xdelta}

For $K=T\times S$ or $G\times S$, we also write $x_i$ for the restriction of
the tautological class to $T^*_\Psi\CB$; on the zero section it is
$-c_1^K(F_i/F_{i-1})$.
The bundle projection $T^*_\Psi\CB\to\CB$ provides a $G\times S$-equivariant deformation retraction
onto the zero section $\CB=G/B$.
Change of groups for the homogeneous space $G/B$ identifies its equivariant cohomology:
\[
 H_{G\times S}^*(T^*_\Psi\CB;\Z)
 \cong H_{B\times S}^*(\pt;\Z)
 \cong\Z[x_1,\dots,x_n,\epsilon],
\]
where the Chern classes $c_j\in\Lambda_{G\times S}$ map to $(-1)^j e_j(x_1,\dots,x_n)$.
For $K=T\times S$ or $G\times S$, let $J_K\subset R_K$ denote the ideal
\begin{align*}
 J_{T\times S}&=\bigl(e_j(x_1,\dots,x_n)-(-1)^je_j(t_1,\dots,t_n)\mid1\leq j\leq n\bigr),\\
 J_{G\times S}&=\bigl(e_j(x_1,\dots,x_n)-(-1)^jc_j\mid1\leq j\leq n\bigr).
\end{align*}
Together with Lemma~\ref{lem:borel}, this gives the canonical presentation
\[
 R_K/J_K\xrightarrow{\sim}H_K^*(T^*_\Psi\CB;\Z),\qquad
 x_i+J_K\longmapsto x_i.
\]
Thus the flag class $x_i$ is represented by the residue class $x_i+J_K$.

\begin{proposition}\label{prop:inclusion_TX_Xdelta_deg_hecke}
For either group $K$, restriction along $\iota_\Psi$ is the surjective
$\Hdeg'$-module homomorphism
\[
 R_K/I_{\rfunc,K}\longrightarrow R_K/J_K,\qquad
 x_i+I_{\rfunc,K}\longmapsto x_i+J_K.
\]
The restricted geometric correspondences are proper over both factors
and act by the operators $\BS_i$ on $R_K/J_K$.
\end{proposition}

\begin{proof}
On the zero section, $q_0(F_i)=F_i$, so the tautological lines restrict to the flag lines.
These classes generate the target over $\Lambda_K$, giving surjectivity
and $I_{\rfunc,K}\subset J_K$. Since $J_K$ is generated by polynomials symmetric in $x_1,\dots,x_n$,
it is preserved by each $\BS_i$, and Theorem~\ref{thm:geometric_hecke}
gives $\Hdeg'$-module compatibility.

Write $X=\SX_\Psi$ and identify $T^*_\Psi\CB$ with its image under $\iota_\Psi$.
For each $s_i\in S'$, every pair $(L_\bullet,L'_\bullet)\in C_i$ has the same terminal lattice $L_n=L'_n$,
and by Proposition~\ref{prop:T^*X_incl} the open condition that $p|_{L_n/L_0}$ is an isomorphism
depends only on $L_n$.
Hence
\[
 \operatorname{pr}_1^{-1}(T^*_\Psi\CB)\cap C_i
 =\operatorname{pr}_2^{-1}(T^*_\Psi\CB)\cap C_i
 =(T^*_\Psi\CB\times T^*_\Psi\CB)\cap C_i,
\]
so the projections $\operatorname{pr}_1|_{C_i},\operatorname{pr}_2|_{C_i}:C_i\to X$ form Cartesian squares over the open immersion $T^*_\Psi\CB\hookrightarrow X$, and their restrictions $C_i\cap(T^*_\Psi\CB\times T^*_\Psi\CB)\to T^*_\Psi\CB$ remain proper.
Proper base change therefore shows that restriction along $\iota_\Psi$ intertwines the action of $\Gamma_i$ with that of its restriction to $T^*_\Psi\CB\times T^*_\Psi\CB$ (when $\epsilon_i-\epsilon_{i+1}\notin\Psi$, the factor $\operatorname{pr}_1^*e_K(\mathcal E_i)$ restricts simultaneously with $[C_i]$); since $\iota_\Psi^*$ is surjective, the restricted correspondences act by $\BS_i$ on $R_K/J_K$.
\end{proof}

\begin{corollary}\label{cor:faithful_hecke}
For either group $K$, the $\Hdeg'$-action on $H_K^*(\SX_{\rfunc};\Z)$ is faithful over $\Z[\epsilon]$.
\end{corollary}

\begin{proof}
For $G\times S$, Proposition~\ref{prop:inclusion_TX_Xdelta_deg_hecke} realizes the faithful polynomial
representation of $\Hdeg'$ on $\Z[x_1,\dots,x_n,\epsilon]$ as a quotient module of $H_{G\times S}^*(\SX_{\rfunc};\Z)$, so Lemma~\ref{lem:polynomial_faithful} proves
the assertion. For $T\times S$, faithfulness follows from the equivariant injection
$H_{G\times S}^*(\SX_{\rfunc};\Z)\hookrightarrow H_{T\times S}^*(\SX_{\rfunc};\Z)$
of Lemma~\ref{lem:borel}.
\end{proof}

\begin{remark}
Faithfulness is taken over $\Z[\epsilon]$ rather than over $\Lambda_K$. For example, on the flag variety $\SX_\varnothing=\CB$, the element $e_1(x_1,\dots,x_n)+c_1$ acts by zero on
$H_{G\times S}^*(\CB;\Z)$, even though it is a nonzero element of $\Lambda_{G\times S}\otimes_{\Z[\epsilon]}\Hdeg'$.
\end{remark}

\begin{corollary}\label{thm:faithful_convolution_embedding}
For either group $K$ and $X=\SX_\Psi$, there is an injective homomorphism of graded $\Z[\epsilon]$-algebras
\[
 \Hdeg'\hookrightarrow\mathcal C_X^\bullet,\qquad
 x_j\longmapsto\Delta_{X*}x_j,\quad
 s_i\longmapsto\Gamma_i\quad(s_i\in S'),\quad
 \epsilon\longmapsto\epsilon[\Delta_X],
\]
where $\Delta_X:X\hookrightarrow X\times X$ is the diagonal embedding and $\Gamma_i$ are the correspondence classes of Section~\ref{sec:geometric_hecke}.
\end{corollary}

\begin{proof}
Under \eqref{eq:ambient_endomorphisms}, these classes act by the polynomial generators
and $\mathcal S_i$. Theorem~\ref{thm:geometric_hecke} and
Corollary~\ref{cor:faithful_hecke} therefore give the homomorphism and its injectivity.
\end{proof}

\subsection{Comparison with Steinberg correspondences}\label{sec:steinberg}

Assume now that $\Psi=\Delta^+$, and set $X=\SX_{\Delta^+}$, $M=T^*\CB$,
and $d=\dim_\C X=n(n-1)$.
After the lattice shift in Remark~\ref{rem:lattice_normalization},
$X$ is the convolution variety $Y(1^n)$ of \cite[Lemma~6.1]{cautis2018quantum}.
Write $Z=M\times_{\CN}M$ for the Steinberg variety and set
\[
 \overline Z=\{(L_\bullet,L'_\bullet)\in X\times X\mid L_n=L'_n\}
\]
with its reduced subscheme structure.
The open embedding of Proposition~\ref{prop:T^*X_incl} identifies
$Z$ with $\overline Z\cap(M\times M)$, as in \cite[Lemma~11.2]{cautis2018quantum}.

For $1\leq i<n$, let $\rho_i:\CB\to\CB_i$ forget the $i$-th subspace and put
$M_i=\rho_i^*T^*\CB_i\subset M$.
The classical simple Steinberg component is
\[
 Z_i=M_i\times_{T^*\CB_i}M_i
\]
(cf.\ \cite[\S12]{ginzburg1998geometric}).
Let $Y_i$ and $\pi_i:Y_i\to B_i$ be the divisor and forgetful morphism
constructed in Section~\ref{sec:geometric_hecke}.

\begin{proposition}\label{prop:compact_steinberg_hecke}
The closure of $Z_i$ in $X\times X$ is
\begin{equation}\label{eq:compact_simple_component}
 \overline Z_i=Y_i\times_{B_i}Y_i.
\end{equation}
It is smooth and irreducible of dimension $d$.
\end{proposition}

\begin{proof}
The variety $Y_i=\SX_{\rfunc_i^-}$ is an iterated projective bundle over a point,
and $Y_i\times_{B_i}Y_i$ is a $\PP^1$-bundle over $Y_i$.
It is therefore smooth and irreducible of dimension $d$.
It lies in $\overline Z$, since both flags have the same terminal lattice.

On $M$, the condition $tL_{i+1}\subset L_{i-1}$ becomes
$\xi F_{i+1}\subset F_{i-1}$, so $Y_i\cap M=M_i$.
By the inverse chart of Proposition~\ref{prop:T^*X_incl},
equality of the retained lattices is equivalent to equality of $\xi$ and the retained flag subspaces.
Thus the intersection of $Y_i\times_{B_i}Y_i$ with $M\times M$ is $Z_i$.
It is a nonempty open subset (containing pairs of equal flags with $\xi=0$),
hence dense in the irreducible fiber product, proving the closure assertion.
\end{proof}

Consequently, the generator images in Corollary~\ref{thm:faithful_convolution_embedding}
are $\Gamma_i=[\Delta_X]+[\overline Z_i]$.
On $M$, the class $\Gamma_i$ restricts to $[\Delta_M]+[Z_i]$,
whose action is identified with $\BS_i$ in Proposition~\ref{prop:inclusion_TX_Xdelta_deg_hecke}.
The corresponding $K$-theoretic construction and open comparison appear in
\cite[Lemma~9.3 and Proposition~11.3]{cautis2018quantum}.

\section{Examples}\label{sec:examples}

We illustrate the divisor construction in rank two and the global projective line
bundle construction for a proper root ideal in rank three.

\subsection{The divisor construction in rank two}

Let $X=\SX_{(0,1)}$ be the compactification of $T^*\PP^1$, and set $\Lambda_{G\times S}=\Z[c_1,c_2,\epsilon]$. By Corollary~\ref{cor:equiv_coho_xdelta_G}, the equivariant cohomology ring is presented as
\[
 H^*_{G\times S}(X;\Z)=\Lambda_{G\times S}[x_1,x_2]/(F_1,F_2),
\]
where
\[
 F_1=x_1^2+c_1x_1+c_2,\qquad
 F_2=x_2^2+c_1x_2+c_2-\epsilon(x_1+x_2+c_1),
\]
with the $\Lambda_{G\times S}$-basis $\{1,x_1,x_2,x_1x_2\}$.
Set $\ell=x_1+x_2+c_1$ and $E=x_2-x_1-\epsilon$; then
$F_2-F_1=E\ell$.
For the smooth divisor $j:Y=\SX_{(0,0)}\cong\PP^1\hookrightarrow X$, we have
\[
 H^*_{G\times S}(Y;\Z)=\Lambda_{G\times S}[x_1,x_2]/(F_1,\ell),\qquad
 j_*[g]=[Eg].
\]
Notice that the divided difference operator $\partial_1$ does not descend to $H^*_{G\times S}(X;\Z)$, because
$\partial_1(F_1)=\ell=x_1+x_2+c_1$ is nonzero in $H^*_{G\times S}(X;\Z)$.
After restriction to the divisor $Y$, however, the relation $\ell=0$ holds, so $\partial_1(F_1)\in(F_1,\ell)$ and the push--pull operator $\pi^*\pi_*$ along
$\pi:Y\to\pt$ acts on $H^*_{G\times S}(Y;\Z)$ as $\partial_1$.
Composing with the Gysin pushforward $j_*$ and adding the identity therefore gives, for every $g\in \Lambda_{G\times S}[x_1,x_2]$,
\[
 (\id+j_*\pi^*\pi_*j^*)[g]=[(\id+E\partial_1)g]=[\BS_1g]\qquad\text{in }H^*_{G\times S}(X;\Z).
\]

\subsection{A projective line bundle for a proper root ideal in rank three}

Let $n=3$ and $\rfunc=(0,1,1)$, so that
$\Psi=\{\epsilon_1-\epsilon_2,\epsilon_1-\epsilon_3\}$.
The simple reflection $s_2$ preserves $\Psi\cup(-\Psi)$, while $\alpha_2=\epsilon_2-\epsilon_3\notin\Psi$.
Set $X=\SX_{(0,1,1)}$, $B_2=\SX_{(0,1,1),2}^P$, $\Lambda_{G\times S}=\Z[c_1,c_2,c_3,\epsilon]$, $f(z)=z^3+c_1z^2+c_2z+c_3$, and
\[
 g(u,v)=\frac{f(u)-f(v)}{u-v}
       =u^2+uv+v^2+c_1(u+v)+c_2.
\]
Then $H^*_{G\times S}(X;\Z)\cong\Lambda_{G\times S}[x_1,x_2,x_3]/(D_1,D_2,D_3)$, where
\[
 \begin{aligned}
 D_1&=f(x_1),\\
 D_2&=f(x_2)-\epsilon g(x_1,x_2),\\
 D_3&=g(x_2,x_3)-\epsilon(x_1+x_2+x_3+c_1).
 \end{aligned}
\]
Because the condition $tL_3\subset L_1$ imposes no constraint on $L_2$, forgetting $L_2$
defines a global $\PP^1$-bundle $\pi_2:X\to B_2$.
The relation ideal $(D_1,D_2,D_3)$ is stable under $\partial_2$, since
\[
 \partial_2D_1=0,\qquad \partial_2D_2=D_3,\qquad \partial_2D_3=0,
\]
and stability follows from the twisted Leibniz rule.
The push--pull operator $\pi_2^*\pi_{2*}$ acts on $H^*_{G\times S}(X;\Z)$ by $\partial_2$, and the twisted relative cotangent line bundle has Euler class
$e_{G\times S}(\Omega^1_{X/B_2}\otimes\C_{-1})=x_3-x_2-\epsilon$.
Thus the second construction of Theorem~\ref{thm:geometric_hecke} gives
\[
 \id+e_{G\times S}(\Omega^1_{X/B_2}\otimes\C_{-1})\cup\pi_2^*\pi_{2*}
 =\id+(x_3-x_2-\epsilon)\partial_2=\BS_2.
\]

\subsection*{Acknowledgments}
The author is deeply grateful to Professor Syu Kato for his guidance, encouragement, and invaluable discussions.

\subsection*{Disclosure of AI use}
The author used GPT-6.1~Sol, GPT-6.2~Sol, and Gemini~4~Argon
to help simplify proofs, check mathematical details, and improve
the English exposition. The author takes full responsibility
for the content of this paper.

\printbibliography
\end{document}